\documentclass[11pt,draft]{amsart}

\usepackage{amsmath}
\usepackage{amssymb}
\usepackage{enumerate}
\usepackage{euscript}
\usepackage{amscd}
\usepackage[all]{xy}

\newtheorem{thm}{Theorem}[section]
\newtheorem{lem}[thm]{Lemma}
\newtheorem{prop}[thm]{Proposition}
\newtheorem{cor}[thm]{Corollary}
\theoremstyle{definition}
\newtheorem{defn}[thm]{Definition}
\newtheorem{rem}[thm]{Remark}

\newtheorem{cond}[thm]{Condition}

\numberwithin{equation}{section}
\theoremstyle{remark}

\newcommand{\bba}{{\mathbb A}}

\newcommand{\bbz}{{\mathbb Z}}

\newcommand{\gA}{{\mathfrak A}}

\newcommand{\gS}{{\mathfrak S}}

\font\tenscr=rsfs10 

\newcommand{\sS}{\hbox{\tenscr S}}

\newcommand{\esB}{{\EuScript{B}}}

\newcommand{\esR}{{\EuScript{R}}}
\newcommand{\esL}{{\EuScript{L}}}

\newcommand{\cA}{{\mathcal A}}
\newcommand{\cE}{{\mathcal E}}

\newcommand{\cN}{{\mathcal N}}

\newcommand{\cT}{{\mathcal T}}

\newcommand{\co}{{\mathcal O}}

\newcommand{\cP}{{\mathcal P}}

\newcommand{\spec}{{\operatorname{Spec}}\,}

\newcommand{\aut}{{\operatorname{Aut}}\,}

\newcommand{\h}{{\operatorname {H}}}

\newcommand{\re}{{\operatorname {Re}}}
\newcommand{\gal}{{\operatorname{Gal}}}
\newcommand{\gl}{{\operatorname{GL}}}

\newcommand{\sst}{{\operatorname{ss}}}
\newcommand{\pfaff}{{\operatorname{Pfaff}}}

\newcommand{\sym}{{\operatorname{Sym}}}

\newcommand{\br}{{\operatorname{B}}}

\newcommand{\A}{\bba}
\newcommand{\Z}{\bbz}

\newcommand{\mk}{k^{\times}}

\newcommand{\ti}{\widetilde}

\newcommand{\twtw}[4]
{\begin{pmatrix}{#1}&{#2}\\{#3}&{#4}\\\end{pmatrix}}
\newcommand{\stwtw}[4]
{\scriptsize\begin{pmatrix}{#1}&{#2}\\{#3}&{#4}\\\end{pmatrix}}

\newcommand{\thth}[9]
{\begin{pmatrix}{#1}&{#2}&{#3}\\{#4}&{#5}&{#6}\\{#7}&{#8}&{#9}\\\end{pmatrix}}

\newcommand{\rh}{{\rm H}}
\newcommand{\ks}{k^{\rm s}}
\newcommand{\sfd}{\mathsf d}

\begin{document}

\title[Parameterizations of rational orbits]
{On parameterizations of rational orbits of some forms of
prehomogeneous vector spaces}

\author[Takashi Taniguchi]{Takashi Taniguchi}
\address{Department of Mathematical Sciences, University of Tokyo}
\email{tani@ms.u-tokyo.ac.jp}
\thanks{Financial support is provided by
Research Fellowships for Young Scientists
of Japan Society for the Promotion of Science.}
\date{\today}

\begin{abstract}
In this paper, we consider the structure of rational orbit
of some prehomogeneous vector spaces originally motivated by
Wright and Yukie.
After we give a refinement of Wright-Yukie's construction,
we apply it to the inner forms of the $D_5$ and $E_7$ types.
Parameterizations of reducible algebras for the cubic cases
are included.
\end{abstract}

\maketitle

\section{Introduction}\label{sec:introduction}

The aim of paper is to give a complement of the series of the work
\cite{wryu, kayu, yukiel} by Yukie, Wright and Kable
on parameterizing separable algebras of degree $2$, $3$, $4$ and $5$
by means of non-degenerate rational orbits
of certain prehomogeneous vector spaces.
The prehomogeneous vector space we treat mainly in this paper
are inner forms of what are called $D_5$ and $E_7$ type.
We also consider the parameterizations of ``reducible algebras''
for the cubic cases ($G_2$, $E_6$ and $E_7$ type.)

Let $k$ be a field of any characteristic
and $\esB$ a quaternion algebra over $k$ with
the canonical involution $\ast$.
We denote by $\cN$ the reduced norm.
Let ${\mathrm M}_n(\esB)$ be the set of $n$ by $n$ matrices
with entries in $\esB$.
This algebra has an involution
$x=(x_{ij})\mapsto x^\iota=(x_{ji}^\ast)$.
Let $\rh_n(\esB)=\{x\in \mathrm M_n(\esB)\mid x^\iota=x\}$.
Then $\gl_n(\esB)$ acts on this space via
$(g,x)\mapsto g\cdot x=gxg^\iota$
for $g\in\gl_n(\esB),\ x\in\rh_n(\esB)$.
In this paper we consider the representation
\begin{equation}\label{eq:pv}
G=\gl_n(\esB)\times \gl_2(k),
\quad
V=\rh_n(\esB)\otimes k^2,
\qquad (n=2,3.)
\end{equation}
We regard this as a representation of an algebraic group $G$ defined over $k$.
This is an inner form of the prehomogeneous vector space
$(\gl_{2n}\times\gl_2,\wedge^{2}k^{2n}\otimes k^2)$
and in the case $\esB\cong{\mathrm M}_2(k)$ they are equivalent.
These are referred to as the $D_5$ type (resp.\ $E_7$ type)
if $n=2$ (resp.\ $n=3$.)
There is a non-zero polynomial $P$ in $V$ such that $P(gx)=\chi(g)P(x)$
for some non-trivial character $\chi$ on $G$.
We put $V^\sst=\{x\in V\mid P(x)\neq0\}$.
Let $\cA_n$ be the set of isomorphism classes of
separable $k$-algebras of degree $n$.
In Section \ref{sec:rationalorbitsD5E7} we prove the following.
\begin{thm}[Theorem \ref{thm:D5_ROD}, Theorem \ref{thm:E7_ROD}]
\label{thm:introthm}
For $n=2,3$, there exist the canonical surjective map
$\Psi\colon G(k)\backslash V^\sst(k)\ni x\mapsto L_x\in\cA_n$
such that the followings hold.
\begin{enumerate}[{\rm (1)}]
\item For $n=2$, $\Psi$ is bijective.
Also for $x\in V^\sst(k)$, the identity component of the stabilizer
$G_x^\circ$ is isomorphic to
$(\esB\otimes L_x)^\times$ as an algebraic group over $k$.
\item For $n=3$, the fiber of $L\in\cA_3$ is parameterized by the set
\[
\left(\mk\cdot{\rm Im}
	(\cN\colon (\esB\otimes L)^\times\rightarrow L^\times)\right)
\backslash L^\times/\aut_k(L).
\]
\end{enumerate}
\end{thm}
For $n=3$, we also determine the structure of $G_x^\circ$
in Theorem \ref{thm:E7_ROD}
and the explicit correspondence of the fibers of $\Psi$
in Theorem \ref{thm:E7_orbit_parameterization}.
We note that \eqref{eq:pv} and representations
in \cite{kayu} numbered as (1) and (3)
have something in common,
and results similar to Theorem \ref{thm:introthm}
have been obtained for those cases.
After Wright-Yukie \cite{wryu} considered $8$ prehomogeneous
vector spaces, some of their $k$-forms were studied
by Kable-Yukie \cite{kayu} and by the author \cite{zfps},
and this paper almost covers the remaining $k$-forms of them.
For the classification of $k$-forms of prehomogeneous vector spaces,
see H. Saito's paper \cite{hsaitoa}.
We briefly discuss the connection of our results
and H. Saito's decomposition formula \cite{hsaitob} of the
global zeta functions of prehomogeneous vector spaces
at the end of Section \ref{sec:rationalorbitsD5E7}.

We follow the argument of \cite{wryu, kayu}, but we give a refinement.
Let $(G,V)$ be a prehomogeneous vector space
treated in \cite{wryu} or \cite{kayu} parameterizing $\cA_i$.
They constructed the map $G(k)\backslash V^\sst(k)\rightarrow \cA_i$
in two ways, via geometry and cohomology and then showed
that two constructions coincide by means of case by case computation.
In this paper we interpret the cohomology theoretic construction
more intrinsically so that it naturally coincides with the
geometric construction.
We see this in Theorem \ref{thm:cohomology_geometry}.
This also simplify the study of the fiber structure of $\Psi$.
We write the proof for the representation \eqref{eq:pv}
but the arguments are parallel for cases
those treated in \cite{wryu, kayu, zfps}.
For this see Remark \ref{rem:8cases}.

Another topic we consider in this paper is on parameterizations
of reducible algebras. Let $\cE_n\subset \cA_n$
be the subset consisting of degree $n$ field extensions of $k$.
From number theoretical interests, we often would like to understand
the set $\cE_n$ rather than $\cA_n$.
For that purpose, it is sometimes useful to parameterize
the set of reducible algebras $\cA_n\setminus \cE_n$
via linear representations.
For the case of the space of binary cubic forms ($G_2$ type)
this was done by Shintani \cite{shintanib}.
In Section 4 we will consider for other cubic cases
i.e., $E_6$ and $E_7$ types.
We hope this result is useful for studying the theory of zeta functions
and some other directions.

This paper is organized as follows.
In Section \ref{sec:galoiscohomology}
we review the definition and basic properties
of the non-abelian Galois cohomology.
In Section \ref{sec:rationalorbitsD5E7} we first
give the compatibility of geometric parameterization
and cohomological parameterization.
We refine the argument of \cite{wryu}
and simplify the description of the fiber structure.
After that we treat $D_5$ and $E_7$ type
prehomogeneous vector spaces explicitly.
A discussion of a decomposition of the global zeta function
is given at the end of the section.
In Section \ref{sec:reducible} we consider
the parameterizations of reducible algebras
for the cubic cases.

\bigskip
{\bf Notation.}
Throughout this paper $k$ is a field of any characteristic.
We fix a separable closure $\ks$ of $k$ and regard any field extensions
of $k$ as subfields of $\ks$.
Let $\Gamma=\gal(\ks/k)$. We assume $\Gamma$ acts on $\ks$ from the right.
Let $X$ be a scheme over $k$.
For a $k$-algebra $R$, the set of $R$-rational points is denoted by $X(R)$
and the base change of $X$ to a scheme over $R$ is denoted by $X\times_kR$
or simply $X_R$. We write $X^{\rm s}=X\times_k\ks$.

For a finite separable field extension $L/k$ and a scheme $Y$ over $L$,
we denote by $\esR_{L/k}(Y)$ the restriction of scalars.
Let $L=L_1\times\dots\times L_n$ where each $L_i/k$
is a finite separable field extension and $X$ a scheme over $k$.
We define a scheme $\esR_{L/k}(X_L)$ over $k$ by
\[
\esR_{L/k}(X_L)
=\esR_{L_1/k}(X_{L_1})\times\dots\times\esR_{L_n/k}(X_{L_n})
\]
This scheme represents the functor $R\mapsto X(L\otimes R)$
where $R$ is any $k$-algebra.

\section{Review of Galois cohomology}\label{sec:galoiscohomology}
In this section we recall the definition of
Galois cohomology and review its properties we use in this paper.
Most of the proofs are omitted. For details, Serre's book
\cite[Chapter I \S5, Chapter III \S1]{serrea}.
We first recall the definition of the Galois cohomology
and properties of sequences.
\begin{defn}
Let $G$ be a group scheme over $k$.
We consider the discrete topology on $G(\ks)$.
A continuous function
$h\colon \Gamma \rightarrow G(\ks)$
is called a {\em 1-cocycle} if
$h(\sigma\tau)=h(\tau)h(\sigma)^\tau$ for all $\sigma,\tau\in\Gamma$.
Two 1-cocycles $h_1$ and $h_2$ are equivalent if there exists $g\in G(\ks)$
such that $h_2(\sigma)=g^{-1}h_1(\sigma)g^\sigma$ for all $\sigma$.
The set of equivalence classes of 1-cocycles is denoted by $\rh^1(k,G)$
and called the {\em first Galois cohomology set} of $G$.
The class of the trivial 1-cocycle
$\Gamma \ni \sigma\mapsto 1\in G(\ks)$
(where $1$ is the identity)
is called the {\em trivial class}.
\end{defn}

By definition, $\rh^1(k,G_1\times G_2)=\rh^1(k,G_1)\times\rh^1(k,G_2)$.
When $G$ is not commutative, $\rh^1(k,G)$ does not have a group structure
in general.
But by regarding the trivial class as its base point,
$\rh^1(k,G)$ becomes a pointed set.
If $G_1\rightarrow G_2$ is a homomorphism of group schemes over $k$,
it obviously induces a homomorphism $\rh^1(k,G_1)\rightarrow \rh^1(k,G_2)$
of pointed sets.
For the rest of this paper we may use the notation
$h=(h_\sigma)_{\sigma\in\Gamma}$ for a 1-cocycle,
where $h_\sigma$ stands for the value of $h$ at $\sigma\in\Gamma$.

Let $G_1,G_2$ and $G_3$ be algebraic groups over $k$. We say a sequence
\begin{equation}\label{eq:sequence}
1\longrightarrow G_1
\overset{\phi_1}{\longrightarrow} G_2
\overset{\phi_2}{\longrightarrow} G_3
\longrightarrow 1
\end{equation}
is {\em exact}
if $G_1$ is a normal subgroup of $G_2$ by $\phi_1$, $\ker(\phi_2)=G_1$,
and $G_2(\ks)\rightarrow G_3(\ks)$ is surjective.
In this situation $G_3(k)$ acts on $\rh^1(k,G_1)$ from the right.
We recall its definition.
Suppose $g\in G_3(k)$ and $\overline h\in\rh^1(k,G_1)$.
Let $h=(h_\sigma)$ be a 1-cocycle in the class $\overline h$.
Take $\tilde g\in G_2(\ks)$ such that $\phi_2(\tilde g)=g$.
Then we see that
$h'_\sigma=\tilde g^{-1}h_\sigma \tilde g^\sigma$
is in $G_1(\ks)$ for all $\sigma\in\Gamma$
and $h'=(h'_\sigma)$ is a 1-cocycle of $G_1$.
It is easy to see that the class $\overline{h'}$
 of $h'$ does not depend on the choice
of $\tilde g$ and $h$, and $(g,\overline h)\mapsto \overline {h'}$
defines a right action of $G_3(k)$.
By considering the action of $G_3(k)$ on the trivial element
in $\rh^1(k,G_1)$, we get a map $G_3(k)\rightarrow \rh^1(k,G_1)$.
\begin{thm}\label{thm:cohomology_exact}
If \eqref{eq:sequence} is exact, then the sequence
\begin{equation*}
1	\longrightarrow G_1(k)
	\longrightarrow G_2(k)
	\longrightarrow G_3(k)
	\longrightarrow \rh^1(k,G_1)
	\longrightarrow \rh^1(k,G_2)
	\longrightarrow \rh^1(k,G_3)
\end{equation*}
is also exact. Moreover, the map
$\rh^1(k,G_1)/G_3(k)\rightarrow \rh^1(k,G_2)$ is injective.
\end{thm}
Note that the exactness in the theorem above means that
the classes which map to the trivial class come from the previous set.

Let $G$ be an affine algebraic group and $H$ a closed subgroup.
Then it is known that the geometric quotient $X=G/H$ exists.
The scheme $X$ does not have a group structure in general
but in the case the injection $G(\ks)/H(\ks)\rightarrow X(\ks)$
is surjective also,
we can still define the map $X(k)\rightarrow \rh^1(k,H)$
the same way as above.
That is, for $x\in X(k)$, take $g\in G(\ks)$
that goes to $x$ by the quotient map 
and then we map $x$ to a cohomology class of the 1-cocycle
$(g^{-1}g^\sigma)_{\sigma\in\Gamma}$.
For this map the following holds.
\begin{thm}\label{thm:rationalorbits_exact}
Assume $G(\ks)/H(\ks)\rightarrow X(\ks)$ is bijective.
Then
$G(k)\backslash X(k)
\rightarrow \rh^1(k,H)$ is injective and
\begin{equation*}
G(k)\backslash X(k)
\longrightarrow \rh^1(k,H)
\longrightarrow \rh^1(k,G)
\end{equation*}
is exact.
\end{thm}

We next recall how the Galois cohomology classifies
$k$-forms of a given algebraic object.
Let $X$ be a variety over $k$. (We do not assume $X$ is connected.)
Then $\Gamma$ acts on $X^{\rm s}$ from the right.
For $\sigma\in\Gamma$ we denote the induced action
$X^{\rm s}\rightarrow X^{\rm s}$ by $\sigma_X$;
\begin{equation*}
\xymatrix{
X^{\rm s}\ar[r]^{\ \ \ \sigma_X\ \ \ } \ar[d]
& X^{\rm s} 
\ar[d]\\
\spec\ks
\ar[r]^{(\sigma^{-1})^\ast}
&
\spec\ks
}
\end{equation*}

Let $\aut(X)$ be the automorphism group of $X$.
We assume $\aut(X)$ acts on $X$ from the right.
It is known that $\aut(X)$ is a group scheme over $k$
and its $R$-rational points $\aut(X)(R)$
is the group of $R$-automorphism of $X\times_kR$.
Hence $\Gamma$ acts on $\aut(X)(\ks)=\aut_{\ks}(X^{\rm s})$ from the right.
For $\alpha\in\aut_{\ks}(X^{\rm s})$ and $\sigma\in\Gamma$,
this action is given by
\[
\alpha^\sigma=\sigma_X\circ\alpha\circ \sigma_X^{-1};\quad
X^{\rm s}	\overset{\sigma_X^{-1}}\longrightarrow
X^{\rm s}	\overset{\alpha}\longrightarrow
X^{\rm s}	\overset{\sigma_X}\longrightarrow
X^{\rm s}.	
\]
\begin{defn}
A variety $Y$ over $k$ is called a {\em $k$-form} of $X$
if $Y^{\rm s}$ is $\ks$-isomorphic to $X^{\rm s}$.
\end{defn}
Let $Y$ be a $k$-form of $X$.
Take a $\ks$-isomorphism $\psi\colon Y^{\rm s}\rightarrow X^{\rm s}$.
For $\sigma\in\Gamma$, let
\[
h_Y(\sigma)=\sigma_X\circ\psi\circ\sigma_Y^{-1}\circ\psi^{-1};\quad
X^{\rm s}	\overset{\psi^{-1}}\longrightarrow
Y^{\rm s}	\overset{\sigma_Y^{-1}}\longrightarrow
Y^{\rm s}	\overset{\psi}\longrightarrow
X^{\rm s}	\overset{\sigma_X}\longrightarrow
X^{\rm s},	
\]
which is an element of $\aut(X)(\ks)=\aut_{\ks}(X^{\rm s})$.
Then since
\begin{align*}
h_Y(\tau)h_Y(\sigma)^\tau
&=	h_Y(\sigma)^\tau\circ h_Y(\tau)
=	\tau_X\circ h_Y(\sigma)\circ\tau_X^{-1}\circ h_Y(\tau)\\
&=	\tau_X\circ \sigma_X\circ \psi\circ \sigma_Y^{-1}\circ\psi^{-1}
	\circ\tau_X^{-1}\circ\tau_X\circ\psi\circ\tau_Y^{-1}\circ\psi^{-1}\\
&=	(\sigma\tau)_X\circ \psi\circ (\sigma\tau)_Y^{-1}\circ \psi^{-1}\\
&=	h_Y(\sigma\tau),
\end{align*}
$h_Y$ is a 1-cocycle of $\aut(X)$. It is easy to see that
its cohomology class does not depend on the choice of $\psi$.
We also denote the cohomology class by $h_Y$.
The following lemma is easy to prove.
\begin{lem}
By associating $h_Y$ to a $k$-form $Y$ of $X$, we have an injective map
from the set of isomorphism classes of $k$-forms of $X$
to $\rh^1(k,\aut(X))$.
In this correspondence, the isomorphism class of $X$ goes to the trivial class.
\end{lem}
In Section \ref{sec:rationalorbitsD5E7} we use the following theorem.
\begin{thm}\label{thm:kform_vs_cohomology}
If $X$ is quasi-projective as a variety (not necessary connected),
then the map above is bijective.
\end{thm}

The following basic result is useful
for the concrete calculation in Section \ref{sec:rationalorbitsD5E7}.
\begin{prop}\label{prop:cohomology_compute}
\begin{enumerate}[{\rm (1)}]
\item
Let $A$ be a finite dimensional simple algebra over $k$.
We regard $A^\times$ as an algebraic group over $k$.
Then $\rh^1(k,A^\times)=\{1\}$.
\item
Let $L/k$ be a finite separable field extension and
$G$ a connected algebraic group over $L$.
Then $\rh^1(k,\esR_{L/k}(G))\cong \rh^1(L,G)$.
\end{enumerate}
\end{prop}

\section{Rational orbits of inner forms of $D_5$, $E_7$ types}
\label{sec:rationalorbitsD5E7}

\subsection{A general construction}
Let $\esB$ be a quaternion algebra over $k$.
The canonical involution of $\esB$ is denoted by $\ast$.
For a $k$-algebra $R$, we put $\esB_R=B\otimes R$.
Let $n$ be a positive integer.
Let ${\mathrm M}_n(\esB)$ be the set of $i$ by $i$ matrices
with entries in $\esB$, which is an associative $k$-algebra.
We denote by $\cT$ and $\cN$ the reduced trace and the reduced norm of
the simple algebra ${\mathrm M}_n(\esB)$.
It is known that $x\in{\mathrm M}_n(\esB)$ is invertible
if and only if $\cN(x)$ is invertible.
The group of invertible elements is denoted by $\gl_n(\esB)$.
For $x=(x_{ij})\in{\mathrm M}_n(\esB)$,
we define $x^\iota=(x_{ji}^\ast)\in{\mathrm M}_n(\esB)$.
Then $\iota$ is an involution of the algebra ${\mathrm M}_n(\esB)$.
Let $\rh_n(\esB)=\{x\in \mathrm M_n(\esB)\mid x^\iota=x\}$,
the space of Hermitian forms of degree $n$.
Then $\gl_n(\esB)$ acts on this space via
\[
(g,x)\longmapsto g\cdot x=gxg^\iota,\qquad
g\in\gl_n(\esB),\ x\in\rh_n(\esB),
\]
which is a linear representation of $\gl_n(\esB)$.
It is known that there is a polynomial of degree $n$
on $\rh_n(\esB)$ which is a square root of $\cN$.
This is uniquely determined
if we assume the value at $1_n\in \rh_n(\esB)$ is $1$.
We denote this polynomial by $\pfaff_n$,
since this coincides with the classical Pfaffian when
$\esB\cong {\mathrm M}_2(k)$.
By definition,
\[
\pfaff_n(1_n)=1,\quad
\pfaff_n(x)^2=\cN(x),\quad
\pfaff_n(g\cdot x)=\cN(g)\pfaff_n(x)
\]
for $g\in\gl_n(\esB),\ x\in\rh_n(\esB)$.

What we are in interest is the representation
\begin{equation}\label{eq:pair_hermitian_repn}
G=\gl_n(\esB)\times \gl_2,
\quad
V=\rh_n(\esB)\otimes k^2
\qquad
(n=2,3),
\end{equation}
where we consider the outer tensor product representation.
The representation $(G,V)$ becomes a prehomogeneous vector space
if and only if $n=1,2,3$. Hence {\em we consider the case $n=2,3$}
for the rest of this paper
(the case $n=1$ is not interesting from our view in this paper.)
We regard this as a representation of algebraic group over $k$.
This is an inner form of
\begin{equation}\label{eq:split_repn}
(\gl_{2n}\times\gl_2,\wedge^{2}k^{2n}\otimes k^2)
\end{equation}
and in the case $\esB\cong{\mathrm M}_2(k)$ they are equivalent.
We say the representation \eqref{eq:pair_hermitian_repn} {\em split}
if $\esB\cong{\mathrm M}_2(k)$ and {\em non-split} otherwise.
It is well known that $\esB\otimes\ks\cong{\mathrm M}_2(\ks)$.
Hence we have:
\begin{prop}\label{prop:ksep_equiv}
The representation \eqref{eq:pair_hermitian_repn}
is equivalent to \eqref{eq:split_repn} over $\ks$.
\end{prop}

We describe the action more explicitly.
Throughout this paper,
we express elements of
$G$ as $g=(g_1,g_2)$ where $g_1\in\gl_n(\esB)$, $g_2\in\gl_2$ and
$V\cong {\mathrm H}_n(\esB)\oplus{\mathrm H}_n(\esB)$ as $x=(x_1,x_2)$
where $x_i\in {\mathrm H}_n(\esB)$.
We call $V$ the space of pairs of binary (resp.\ ternary) Hermitian forms
if $n=2$ (resp.\ $n=3$).
We identify $x=(x_1,x_2)\in V$ with
the Hermition form $x(v)=x_1v_1+x_2v_2$
of linear form in the variables $v_1$ and $v_2$,
which we collect into the row vector $v=(v_1,v_2)$.
Then the left action of $g\in G$ on $x\in V$ is defined by
\begin{equation*}
(g\cdot x)(v)=g_1\cdot x(vg_2)=g_1x(vg_2)g_1^\iota,
\end{equation*}
or, in other words,
\[
(g_1,\begin{pmatrix} a&b\\c&d\end{pmatrix})\cdot(x_1,x_2)
	=(g_1(ax_1+bx_2)g_1^\iota,g_1(cx_1+dx_2)g_1^\iota).
\]
We put $F_x(v)=\pfaff_n(x(v))$.
This is a binary quadratic form (resp.\ cubic form)
in variables $v=(v_1,v_2)$ if $n=2$ (resp.\ $n=3$),
and the discriminant $P(x)\ (x\in V)$ is a polynomial in $V$.
The polynomial $P(x)$ is characterized by 
\begin{equation*}
P(x)=\prod_{1\leq i<j\leq n}(\alpha_i\beta_j-\alpha_j\beta_i)^2
\quad
\text {for}
\quad
F_x(v)=\prod_{1\leq i\leq n}(\alpha_iv_1-\beta_iv_2), \ \ x\in V(\ks).
\end{equation*}
We define $\chi(g)=\cN(g_1)\det(g_2)$ for $n=2$ and
$\chi(g)=\cN(g_1)^2\det(g_2)^3$ for $n=3$.
Then one can easily see that
\begin{equation*}
P(gx)=\chi(g)^2P(x)
\end{equation*}
and hence $P(x)$ is a relative invariant polynomial
with respect to the character $\chi^2$.
Let $S=\{x\in V\mid P(x)=0\}$ and $V^\sst=\{x\in V\mid P(x)\not=0\}$
and call them the set of unstable points and semi-stable points, respectively.
That is, $x\in V$ is semi-stable if and only if
$F_x(v)$ does not have a multiple root in ${\mathbb P}^1=\{(v_1:v_2)\}$.

We now consider rational points of $V^\sst$ and its stabilizer structure.
\begin{defn}
\begin{enumerate}[{\rm (1)}]
\item
Let $\cP_2$ (resp.\ $\cP_3)$ be the set of equivalence classes of $k$-forms of
the variety $\mathrm{Proj}\, k[v_1,v_2]/(v_1v_2)$
(resp.\ $\mathrm{Proj}\, k[v_1,v_2]/(v_1v_2(v_1-v_2)$).
\item
For $n=2,3$, let $\cA_n$ be set of isomorphism classes
of separable $k$-algebras of dimension $n$.
\end{enumerate}
\end{defn}
We note that the contravariant functor
\[
s\colon \cP_n\longrightarrow \cA_n,\quad
Z\longmapsto s(Z)=\Gamma(Z,\co_Z)
\]
gives a categorical equivalence between $\cP_n$ and $\cA_n$.
\begin{defn}
For $x\in V^\sst(k)$, we define
\begin{align*}
Z_x	&=\mathrm{Proj}\, k[v_1,v_2]/(F_x(v)),\\
L_x	&=\Gamma(Z_x,\co_{Z_x}).
\end{align*}
\end{defn}
Since $V^\sst(k)$ is the set of $x$ such that
$F_x$ does not have a multiple root,
$Z_x\in\cP_n$ is a variety over $k$ and $L_x\in\cA_n$.
By definition $Z_x$ depends only on the $G(k)$-orbit of $x$
and hence we have a map $G(k)\backslash V^\sst(k)\rightarrow \cP_n$.
To describe the fiber structure, Wright and Yukie \cite{wryu}
interpreted this map in terms of Galois cohomology.
We give a refined version of their argument.

We fix an arbitrary element $w\in V^\sst(k)$.
By definition, $\cP_n$ is the set of equivalence classes of $k$-forms of $Z_w$.
Hence by Theorem \ref{thm:kform_vs_cohomology} there is the canonical bijection
\begin{equation}
d_w\colon \rh^1(k,\aut(Z_w))\longrightarrow \cP_n.
\end{equation}

Let $G_w=\{g\in G\mid g\cdot w=w\}$ be the stabilizer
and $G_w^\circ$ its identity component.
If $g=(g_1,g_2)\in G_w$, then $\cN(g_1)F_w(vg_2)=F_w(v)$.
Hence each element of $G_w$ defines an automorphism of $Z_w$ and the map
\begin{equation}\label{eq:c_w}
c_w\colon G_w\longrightarrow \aut(Z_w),\quad
g=(g_1,g_2)\longmapsto (Z_w\ni v\longmapsto vg_2\in Z_w)
\end{equation}
gives a homomorphism of algebraic groups.
We note that we assume $\aut(Z_w)$ acts on $Z_w$ from the right.

Let us consider the following condition:
\begin{cond}\label{cond:pv_structure}
\begin{enumerate}[{\rm (1)}]
\item\label{cond:quotient_isom}
For any $w\in V^\sst(k)$, the map $g\mapsto gw$ gives an isomorphism
$G/G_w\rightarrow V^\sst$ as varieties over $k$.
We denote by $a_w$ the inverse map;
\begin{equation}
a_w\colon V^\sst\longrightarrow G/G_w.
\end{equation}
\item\label{cond:singleorbit}
The canonical injection $G(\ks)/G_w(\ks)\rightarrow (G/G_w)(\ks)$
is bijective.
\item\label{cond:stab_exact}
For any $w\in V^\sst(k)$, the sequence
\[
1	\longrightarrow G_w^\circ
	\longrightarrow G_w
	\overset{c_w}{\longrightarrow} \aut(Z_w)
	\longrightarrow 1
\]
is exact.
\item\label{cond:G_coh_triv}
We have $\rh^1(k,G)=\{1\}$.
\end{enumerate}
\end{cond}

\begin{prop}\label{prop:pv_structure}
For the prehomogeneous vector space \eqref{eq:pair_hermitian_repn},
Condition \ref{cond:pv_structure} hold.
\end{prop}
\begin{proof}
(\ref{cond:quotient_isom}) and (\ref{cond:singleorbit})
are proved by Yukie \cite{yukieq}
for a general irreducible regular prehomogeneous vector space.
For (\ref{cond:stab_exact}) we may assume $k=\ks$.
Then $(G,V)$ is the split form.
Since the exactness is proved in \cite{wryu}
for one element, the same is true for any its $G(k)$-orbit.
By (\ref{cond:singleorbit}) $V^\sst(k)$ is a single $G(k)$-orbit
and we have (\ref{cond:stab_exact}).
(\ref{cond:G_coh_triv})
is a consequence of Proposition \ref{prop:cohomology_compute} (2).
\end{proof}

From (\ref{cond:singleorbit}) and (\ref{cond:G_coh_triv})
of Condition \ref{cond:pv_structure},
we have the following by Theorem \ref{thm:rationalorbits_exact}.
\begin{prop}The canonical map
\begin{equation}
b_w\colon G(k)\backslash (G/G_w)(k)\longrightarrow \rh^1(k,G_w)
\end{equation}
is bijective.
\end{prop}
\begin{defn}
The map $a_w$ induces a bijective map
$G(k)\backslash V^\sst(k)\rightarrow G(k)\backslash (G/G_w)(k)$,
which we also denote by $a_w$.
Also the map $\rh^1(k,G_w)\rightarrow \rh^1(k,\aut(Z_w))$
induced by $c_w$ is also denoted by $c_w$.
We define $\Phi_w$ by the following composition;
\begin{equation}
\Phi_w=d_w\circ c_w\circ b_w\circ a_w\colon\ \ 
G(k)\backslash V^\sst(k) \longrightarrow \cP_n.
\end{equation}
\end{defn}
\begin{thm}\label{thm:cohomology_geometry}
For $x\in V^\sst(k)$, we have $\Phi_w(G(k)\cdot x)=Z_x$.
Especially $\Phi_w$ does not depends on the choice of $w\in V^\sst(k)$.
\end{thm}
\begin{proof}
By (\ref{cond:quotient_isom}) and (\ref{cond:singleorbit})
of Condition \ref{cond:pv_structure},
$V^\sst(\ks)$ is a single $G(\ks)$-orbit. Take
$g\in G(\ks)$ such that
\begin{equation}\label{eq:x=gw}
x=g\cdot w,\qquad g=(g_1,g_2).
\end{equation}

Then $b_w(a_w(x))\in \rh^1(k,G_w)$ is the class of the 1-cocycle
$(g^{-1}g^\sigma)_{\sigma\in\Gamma}$.
Hence $c_w(b_w(a_w(y)))\in \rh^1(k,\aut(Z_w))$ 
is the class of the 1-cocycle $(c_w(g^{-1}g^\sigma))_{\sigma\in\Gamma}$.

On the other side since $x=g\cdot w$,
$F_x(v)$ coincides with $F_w(vg_2)$ up to non-zero constant.
Hence
\begin{equation}\label{eq:ks_isom_by_g}
\psi\colon Z_x\times_k \ks \longrightarrow Z_w\times_k\ks,\qquad
v\longmapsto \psi(v)=vg_2
\end{equation}
is an isomorphism of varieties over $\ks$.
Let
\[
h_{Z_x}(\sigma)=
\sigma_{Z_w}\circ \psi\circ \sigma_{Z_x}^{-1}\circ \psi^{-1}
	\in\aut_{\ks}(Z_w\times_k \ks).
\]
Then $d_w^{-1}(Z_x)$ is the cohomology
class of $(h_{Z_x}(\sigma))_{\sigma\in\Gamma}$.
Also for $v\in (Z_w\times_k\ks)\subset {\mathbb P}^1$,
\begin{align*}
h_{Z_x}(\sigma)(v)
&=	\left(\sigma_{Z_w}\circ \psi\circ \sigma_{Z_x}^{-1}\right)(vg_2^{-1})
=	(\sigma_{Z_w}\circ \psi)
	\left(v^{\sigma^{-1}}(g_2^{-1})^{\sigma^{-1}}\right)\\
&=	\sigma_{Z_w}
	\left(v^{\sigma^{-1}}(g_2^{-1})^{\sigma^{-1}}g_2\right)
=	vg_2^{-1}g_2^\sigma,
\end{align*}
which shows that $h_{Z_x}(\sigma)=c_w(g^{-1}g^\sigma)$.
Hence we have $(c_w\circ b_w\circ a_w)(x)=d_w^{-1}(Z_x)$
and this finishes the proof.
\end{proof}
Since $\Phi_w$ does not depend on $w$, we drop the subscript $w$
and write $\Phi$ from now on. As a corollary to the theorem above,
we obtain the following;
\begin{cor}\label{cor:fiber_structure}
For $x\in V^\sst(k)$, the fiber $\Phi^{-1}(Z_x)$
is parameterized by the set:
\[
\rh^1(k,G_x^\circ)/\aut(Z_x)(k),
\]
where the action in the right hand side is determined by
the exact sequence in Condition \ref{cond:pv_structure}
(\ref{cond:stab_exact}).
Moreover, if a cohomology class in $\rh^1(k,G_x^\circ)$
is represented by 1-cocycle $(g^{-1}g^\sigma)_{\sigma\in\Gamma}$
for some $g\in G(\ks)$, then the corresponding element
in $\Phi^{-1}(Z_x)$ is the orbit of $gx\in V^\sst(k)$.
\end{cor}
\begin{proof}
We regard $\Phi$ as the composition
$d_x\circ c_x\circ b_x\circ a_x$.
Since $d_x$, $b_x$ and $a_x$ are bijective,
we consider the fiber of
$c_x\colon \rh^1(k,G_x)\rightarrow \rh^1(k,\aut(Z_x))$.
Since $d_x^{-1}(Z_x)\in \rh^1(k,\aut(Z_x))$ is the trivial class,
the first statement follows from the exactness of the sequence
\[
\rh^1(k,G_x^\circ)
\longrightarrow		\rh^1(k,G_x)
\overset{c_x}{\longrightarrow}		\rh^1(k,\aut(Z_x))
\]
and the injectivity of the map
$\rh^1(k,G_x^\circ)/\aut(Z_x)(k)
\rightarrow		\rh^1(k,G_x)$
those asserted in Theorem \ref{thm:cohomology_exact}.
The second statement follows from
the definitions of maps of the Galois cohomology.
\end{proof}
\begin{rem}\label{rem:8cases}
We prove results above using Condition \ref{cond:pv_structure}
and formulae \eqref{eq:c_w}, \eqref{eq:x=gw}, \eqref{eq:ks_isom_by_g}.
It is written explicitly or implicitly
in \cite{wryu, kayu} that all the prehomogeneous vector
spaces they considered satisfy the same condition
for certain $Z_x$ defined geometrically,
hence our arguments are applicable for their cases.
In \cite{wryu, kayu} they chose a specific element $w$ so that
the exact sequence of \ref{cond:pv_structure} (\ref{cond:stab_exact}) splits.
Our approach does not need to choose such an element, and hence fitting
in cases that there does not exist $w$ such that $c_w$ splits
as we treated in \cite{zfps}.


\end{rem}

\begin{defn}
We put $\Psi=s\circ \Phi\colon G(k)\backslash V^\sst(k)\rightarrow \cA_n$.
\end{defn}

For $k$-algebra $R$, the norm map $\esB_R\rightarrow R$
induced from $\cN\colon \esB\rightarrow k$ is also denoted by $\cN$.
We often regard $\esB^\times$ as an algebraic group over $k$.
We hope the meaning is clear from the context.

\subsection{The case of $D_5$ type}
We now consider the case $n=2$ and $n=3$ separately.
In this subsection we consider the case $n=2$, i.e.,
the space of pairs of binary Hermitian forms.
Let
\[
w=\left(\twtw 0001,\twtw 1000\right)\in V(k).
\]

\begin{prop}\label{prop:D5_neutral_element}
We have $w\in V^\sst(k)$, $L_w=k\times k$ and
\[
G_w^\circ=
\left\{\left(\twtw {b_1}00{b_2},\twtw {t_1}00{t_2}\right)
\ \vrule\ 
\begin{array}{l}
b_i\in \esB^\times, t_i\in \mathbb G_m,\\
t_1\cN(b_1)=t_2\cN(b_2)=1\\\end{array}
\right\}.
\]
\end{prop}
\begin{proof}
Since $F_w(v)=v_1v_2$, we have $w\in V^\sst(k)$ and $L_w=k\times k$.
Let $H$ denote the group of the right hand side.
We immediately see that elements of $H$ stabilize $w$
and that $H$ is connected. Hence $H\subset G_w^\circ$.
Since
\[
\dim H=8=\dim G-\dim V=\dim G-\dim(G/G_w)
=\dim G_w=\dim G_w^\circ
\]
and both $H$ and $G_w^\circ$ are connected, they coincide.
\end{proof}

Let $L=k(\alpha)$ is a quadratic field extension.
Let $\nu$ be the generator of $\gal(L/k)$.
Put $a_1=\alpha+\alpha^\nu$, $a_2=\alpha\alpha^\nu$ and
\begin{equation*}
x_\alpha=\left(	\twtw 0{1}{1}{a_1},
		\twtw {1}{a_1}{a_1}{a_1^2-2a_2}
	\right)
\in V(k).
\end{equation*}
Let
\begin{equation*}
g_\alpha=\left(\twtw 11{\alpha}{\alpha^\nu},
	\frac{1}{\alpha-\alpha^\nu}\twtw {-1}{1}{-\alpha^\nu}{\alpha}
\right)\in G(L).
\end{equation*}
By computation we have $x_\alpha=g_\alpha w$.

\begin{prop}\label{prop:D5_quadratic_element}
We have $x_\alpha\in V^\sst(k)$, $L_{x_\alpha}=L$, and
$G_{x_\alpha}^\circ\cong \esR_{L/k}(\esB_L^\times)$.
\end{prop}
\begin{proof}
Since $F_{x_\alpha}(v)=-(v_1^2+a_1v_1v_2+a_2v_2^2)=
-(v_1+\alpha v_2)(v_1+\alpha^\nu v_2)$,
we have $x_\alpha\in V^\sst(k)$ and $L_{x_\alpha}=L$.
We consider $G_x^\circ$. Let $R$ be any $k$-algebra.
Then since $x=g_\alpha w$, 
$G_x^\circ(R\otimes L)=g_\alpha G_w^\circ(R\otimes L)g_\alpha^{-1}$.
Also $G_x^\circ(R)=
\left\{p\in G_x^\circ(R\otimes L) \mid p^\nu=p\right\}$.
Using $g_\alpha^\nu=g_\alpha\tilde\nu$ where
\[
\tilde\nu=\left(\twtw 0110,\twtw0110\right),
\]
by computation we have
\[
G_x^\circ(R)=g_\alpha
\left\{\left(\twtw {b}00{b^\nu},\twtw {t}00{t^\nu}\right)
\ \vrule\ 
\begin{array}{l}
b\in (\esB_L\otimes R)^\times,\\ t\in (L\otimes R)^\times,
t\cN(b)=1\\\end{array}
\right\}g_\alpha^{-1}.
\]
Hence $G_x^\circ(R)\cong (\esB_L\otimes R)^\times$
and this isomorphism is functorial with respect to $k$-algebras.
This proves the proposition.
\end{proof}

\begin{prop}
The map $\Psi\colon G(k)\backslash V^\sst(k)\rightarrow \cA_2$ is bijective.
\end{prop}
\begin{proof}
By Propositions \ref{prop:D5_neutral_element} and
\ref{prop:D5_quadratic_element}, $\Psi$ is surjective.
By Proposition \ref{prop:cohomology_compute},
$\rh^1(k,\esB^\times\times\esB^\times)
=\rh^1(k,\esR_{L/k}(\esB_L^\times))=\{1\}$
where $L/k$ is a quadratic field extension.
This combined with
Propositions \ref{prop:D5_neutral_element} and
\ref{prop:D5_quadratic_element}
shows that for any $L\in \cA_2$, there exists $x\in V^\sst(k)$
such that $L_x=L$ and $\rh^1(k,G_x^\circ)=\{1\}$.
Hence by Corollary \ref{cor:fiber_structure} each fiber of
elements of $\cA_2$ consists of a single orbit. Hence $\Psi$ is injective.
\end{proof}

We summarize the result of this subsection as follows.
\begin{thm}\label{thm:D5_ROD}
 Let $n=2$.
\begin{enumerate}[{\rm (1)}]
\item
The map
\[
\Psi\colon G(k)\backslash V^\sst(k)\longrightarrow \cA_2,\quad
G(k)\cdot x\longmapsto L_x
\]
is bijective.
\item
For $x\in V^\sst(k)$, $G_x^\circ\cong \esR_{L_x/k}(\esB_{L_x}^\times)$.
\end{enumerate}
\end{thm}

\subsection{The case of $E_7$ type: Statements}
The case $n=3$ becomes a little more complicated than the case $n=2$,
and we give statements and proofs separately.
We give statements in this subsection and prove in the next subsection.
For the space of pairs of ternary Hermitian forms, the following theorems hold.
\begin{thm}\label{thm:E7_ROD}
Let $n=3$.
\begin{enumerate}[{\rm (1)}]
\item
The map
\[
\Psi\colon G(k)\backslash V^\sst(k)\longrightarrow \cA_3,\quad
G(k)\cdot x\longmapsto L_x
\]
is surjective, and the fiber of any $L\in\cA_3$ is parameterized by the set
\[
\left(\mk\cdot{\rm Im}
	(\cN\colon \esB_L^\times\rightarrow L^\times)\right)
\backslash L^\times/\aut_k(L).
\]
\item
Let $x\in V^\sst(k)$.
There exists an exact sequence of algebraic groups
\[
1
\longrightarrow		G_x^\circ
\longrightarrow		{\mathbb G}_m\times \esR_{L_x/k}(\esB_{L_x}^\times)
\overset{\varphi_x}{\longrightarrow}	\esR_{L_x/k}({\mathbb G}_m)
\longrightarrow		1,
\]
where $\varphi_x$ is given by $\varphi_x(t,b)=t\cN(b)$.
\end{enumerate}
\end{thm}
A parameterization in Theorem \ref{thm:E7_ROD} (1) is given by:
\begin{thm}\label{thm:E7_orbit_parameterization}
\begin{enumerate}[{\rm (1)}]
\item
Let $L=k\times k\times k$.
Then $(p_1,p_2,p_3)\in (k^\times)^3=L^\times$ corresponds to the orbit of
\[
\left(\thth 0000{p_2}000{-p_3}, \thth {p_1}000{-p_2}0000\right).
\]
\item
Let $L=k\times F$ where $F=k(\alpha)$ is a quadratic field extension of $k$.
Let $\nu$ be the generator of $\gal(F/k)$.
Then $(p,\lambda)\in k^\times\times F^\times=L^\times$ corresponds
to the orbit of
\[
\left(
\thth 0000{\Lambda_0}{\Lambda_1}0{\Lambda_1}{\Lambda_2},
\thth {-p}000{\Lambda_1}{\Lambda_2}0{\Lambda_2}{\Lambda_3}\right),
\]
where we put
$\Lambda_i=
(\lambda\alpha^i-(\lambda\alpha^i)^\nu)/(\alpha-\alpha^\nu)\in k$.
\item
Let $L=k(\beta)$ be a cubic field extension of $k$.
We denote by $L^{\rm nc}$ the normal closure of $L$,
and $\beta_1=\beta$, $\beta_2=\beta^\tau$, $\beta_3=\beta^\mu$
($\tau,\mu\in\gal(L^{\rm nc}/k)$)
the all conjugates of $\beta$.
Then $\delta\in\L^\times$ corresponds to the orbit of
\[
\left(
\thth {\Delta_0}{\Delta_1}{\Delta_2}
	{\Delta_1}{\Delta_2}{\Delta_3}{\Delta_2}{\Delta_3}{\Delta_4},
\thth {\Delta_1}{\Delta_2}{\Delta_3}
	{\Delta_2}{\Delta_3}{\Delta_4}{\Delta_3}{\Delta_4}{\Delta_5}
\right),
\]
where we put $\delta_1=\delta$, $\delta_2=\delta^\tau$, $\delta_3=\delta^\mu$
and
\[
\Delta_i=
\frac{
	 \delta_1\beta_1^i(\beta_2-\beta_3)
	+\delta_2\beta_2^i(\beta_3-\beta_1)
	+\delta_3\beta_3^i(\beta_1-\beta_2)
	}
{(\beta_1-\beta_2)(\beta_2-\beta_3)(\beta_3-\beta_1)}\in k.
\]
\end{enumerate}
\end{thm}

\subsection{The case of $E_7$ type: Proofs}
We now give proofs of Theorems \ref{thm:E7_ROD}
and \ref{thm:E7_orbit_parameterization}
in order. We keep the notation in the previous subsection.
As for the $D_5$ case,
we first choose a specific element in the fiber for each $L\in\cA_3$
and determine its stabilizer structure.

We write ``diagonal'' elements of $G$ as
\[
\sfd(b_1,b_2,b_3,t)=\left(\thth {b_1}000{b_2}000{b_3},\twtw t00t\right)\in G
\qquad
(b_i\in\esB^\times, t\in\mathbb G_m.)
\]
Let
\begin{equation}\label{eq:E7_w}
w=\left(\thth 00001000{-1},\thth 1000{-1}0000\right)\in V(k).
\end{equation}
We define
\[
K^w=\{\sfd(b_1,b_2,b_3,t)\mid b_i\in\esB^\times, t\in\mathbb G_m\}\subset G
\]
and
\[
\varphi_w\colon K^w\longrightarrow \mathbb G_m^3,\quad
\sfd(b_1,b_2,b_3,t)\longmapsto (t\cN(b_1),t\cN(b_2),t\cN(b_3)).
\]
Then the following can be proved similarly to Proposition
\ref{prop:D5_neutral_element}.
\begin{prop}\label{prop:E7_neutral_element}
We have $w\in V^\sst(k)$, $L_w=k\times k\times k$ and
$G_w^\circ=\ker(\varphi_w)$.
\end{prop}

Let $L=k\times F\in\cA_3$
where $F=k(\alpha)/k$ is a quadratic field extension.
Put $a_1=\alpha+\alpha^\nu$, $a_2=\alpha\alpha^\nu$ and
\begin{equation}\label{eq:E7_xalpha}
x_\alpha=\left(	\thth 00000101{a_1},
		\thth {-1}0001{a_1}0{a_1}{a_1^2-a_2}
	\right)
\in V(k).
\end{equation}
Let
\[
g_\alpha=
\left(
	\thth 1000110\alpha{\alpha^\nu},
	\frac{1}{\alpha-\alpha^\nu}
		\twtw 10{\alpha^\nu}{\alpha^\nu-\alpha}
\right)\in G(F).
\]
By computation we have $x_\alpha=g_\alpha w$.
Let
\[
K^{x_\alpha}(k)=g_\alpha
\{
\sfd(b_1,b_2,b_2^\nu,t)
	\mid b_1\in\esB^\times,\, b_2\in \esB_F^\times,\, t\in\mk
\}
g_\alpha^{-1}.
\]
\begin{lem}
We have
$K^{x_\alpha}(k)\subset G(k)$.
\end{lem}
\begin{proof}
Let $p=g_\alpha\sfd(b_1,b_2,b_2^\nu,t)g_\alpha^{-1}\in K^{x_\alpha}(k)$.
We already know $p\in G(F)$.
Since $g_\alpha^\nu=g_\alpha\ti\nu$ where
\[
\ti\nu
=\left(
	\thth 100001010,
	\twtw {-1}011
\right),
\]
we have $p^\nu=
g_\alpha\ti\nu\sfd(b_1,b_2^\nu,b_2,t)\ti\nu^{-1}g_\alpha^{-1}
=g_\alpha\sfd(b_1,b_2,b_2^\nu,t)g_\alpha^{-1}=p$.
Hence $p\in G(k)$.
\end{proof}

We regard $K^{x_\alpha}(k)$ as $k$-rational points
of the algebraic group $K^{x_\alpha}\subset G$ over $k$ in an obvious manner.
That is, for any $k$-algebra $R$,
\[
K^{x_\alpha}(R)=g_\alpha
\{
\sfd(b_1,b_2,b_2^\nu,t)\mid
	b_1\in\esB_R^\times,\,
	b_2\in \esB_{F\otimes R}^\times,\,
	t\in R^\times
\}
g_\alpha^{-1},
\]
which we can prove is a subgroup of $G(R)$ similarly to the lemma above.
We define
\[
\varphi_{x_\alpha}\colon
	K^{x_\alpha}\longrightarrow \mathbb G_m\times \esR_{F/k}(\mathbb G_m),
\quad
	g_\alpha\sfd(b_1,b_2,b_2^\nu,t)g_\alpha^{-1}
		\longmapsto (t\cN(b_1), t\cN(b_2)).
\]

\begin{prop}\label{prop:E7_quadratic_element}
We have $x_\alpha\in V^\sst(k)$, $L_{x_\alpha}=L$ and
$G_{x_\alpha}^\circ=\ker(\varphi_{x_\alpha})$.
\end{prop}
\begin{proof}
Since
\[
F_{x_\alpha}(v)=v_2(v_1^2+a_1v_1v_2+a_2v_2^2)=v_2(v_1+\alpha v_2)(v_1+\alpha^\nu v_2),
\]
we have $x_\alpha\in V^\sst(k)$, $L_{x_\alpha}=L$.
Let $R$ be a $k$-algebra.
We consider $G_{x_\alpha}^\circ(R)$.
Since $x_\alpha=g_\alpha w$,
$G_{x_\alpha}^\circ(R\otimes F)=g_\alpha G_{w}^\circ(R\otimes F)g_\alpha^{-1}$.
Let
\[
p=g_\alpha\sfd(b_1,b_2,b_3,t)g_\alpha^{-1}\in G_{x_\alpha}^\circ(R\otimes F),
\quad
b_i\in(\esB_F\otimes R)^\times,\ t\in (F\otimes R)^\times,\ t\cN(b_i)=1.
\]
Then $p=p^\nu$ if and only if $b_1=b_1^\nu$, $b_3=b_2^\nu$,
and $t=t^\nu$.
This shows $G_{x_\alpha}^\circ(R)=\ker(\varphi_{x_\alpha})(R)$.
\end{proof}

Let $L=k(\beta)/k$ be a cubic field extension.
Put $b_1=\beta_1+\beta_2+\beta_3$,
$b_2=\beta_2\beta_3+\beta_3\beta_1+\beta_1\beta_2$,
$b_3=\beta_1\beta_2\beta_3$,
$D=(\beta_1-\beta_2)(\beta_2-\beta_3)(\beta_3-\beta_1)$
and
\begin{equation}\label{eq:E7_xbeta}
x_\beta=\left(
	\thth 00101{b_1}1{b_1}{b_1^2-b_2},
	\thth 01{b_1}1{b_1}{b_1^2-b_2}{b_1}{b_1^2-b_2}{b_1^3-2b_1b_2+b_3}
	\right)
\in V(k).
\end{equation}
Let
\[
g_\beta=
\left(
	\thth 111{\beta_1}{\beta_2}{\beta_3}{\beta_1^2}{\beta_2^2}{\beta_3^2},
	\frac{1}{D}
	\twtw	{\beta_2-\beta_1}{\beta_2-\beta_3}
		{\beta_3(\beta_2-\beta_1)}{\beta_1(\beta_2-\beta_3)}
\right)\in G(L^{\rm nc}).
\]
By computation we have $x_\beta=g_\beta w$.
Let
\[
K^{x_\beta}(k)=g_\beta
\{
\sfd(b,b^\tau,b^\mu,t)
	\mid b\in\esB_L^\times,\, t\in\mk
\}
g_\beta^{-1}.
\]
\begin{lem}
We have
$K^{x_\beta}(k)\subset G(k)$.
\end{lem}
\begin{proof}
Let $L/k$ be a non-Galois cubic extension.
By changing $\tau,\mu\in\gal(L^{\rm nc}/k)$ if necessary,
we assume they are of index $2$.
Then since $\tau(\beta_3)=\beta_3$ and $\mu(\beta_2)=\beta_2$,
we have $g_\beta^\tau=g_\beta\ti\tau$ and $g_\beta^\mu=g_\beta\ti\mu$ where
\[
\ti\tau=
\left(\thth 010100001, \twtw 110{-1}\right),
\quad
\ti\mu=
\left(\thth 001010100, \twtw 0{-1}{-1}0\right).
\]
Also $\mu\tau\mu=\tau\mu\tau\in\gal(L^{\rm nc}/k)$ fixes elements of $L$.
Hence we can see any element
$p=g_\beta\sfd(b,b^\tau,b^\mu,t)g_\beta^{-1}\in G(L^{\rm nc})$
is fixed by $\tau,\mu$, and consequently $p\in G(k)$.
The case $L/k$ is Galois is similarly proved and we omit the repetition here.
\end{proof}
As before we regard $K^{x_\beta}(k)$ as $k$-rational points
of the algebraic group $K^{x_\beta}\subset G$ over $k$.
We define
\[
\varphi_{x_\beta}\colon
	K^{x_\beta}\longrightarrow \esR_{L/k}(\mathbb G_m),
\quad
	g_\beta\sfd(b,b^\tau,b^\mu,t)g_\beta^{-1}
		\longmapsto t\cN(b).
\]

\begin{prop}\label{prop:E7_cubic_element}
We have $x_\beta\in V^\sst(k)$, $L_{x_\beta}=L$ and
$G_{x_\beta}^\circ=\ker(\varphi_{x_\beta})$.
\end{prop}
\begin{proof}
Since
\[
F_{x_\beta}(v)=-(v_1^3+b_1v_1^2v_2+b_2v_1v_2^2+b_3v_2^3)
=-(v_1+\beta_1v_2)(v_1+\beta_2v_2)(v_1+\beta_3v_2),
\]
we have $x_\beta\in V^\sst(k)$, $L_{x_\beta}=L$.
We determine $G_{x_\beta}^\circ$.
Consider the case $L/k$ a non-Galois cubic extension.
We use the notation in the lemma above.
Since $x_\beta=g_\beta w$,
$G_{x_\beta}^\circ(L^{\rm nc})=g_\beta G_{w}^\circ(L^{\rm nc})g_\beta^{-1}$.
Let
\[
p=g_\alpha\sfd(b_1,b_2,b_3,t)g_\alpha^{-1}\in G_{x_\alpha}^\circ(L^{\rm nc}),
\quad
b_i\in\esB_{L^{\rm nc}}^\times,\ t\in (L^{\rm nc})^\times,\ t\cN(b_i)=1.
\]
Then $p=p^\tau=p^\mu$ if and only if
$b_2=b_1^\tau$, $b_3=b_1^\mu$, $b_1=b_1^{\mu\tau\mu}=b_1^{\tau\mu\tau}$
and $t=t^\tau=t^\mu$.
This shows $G_{x_\beta}^\circ(k)=\ker(\varphi_{x_\beta})(k)$.
The equality $G_{x_\beta}^\circ(R)=\ker(\varphi_{x_\beta})(R)$
for a $k$-algebra $R$ is proved in the same way and thus we have
$G_{x_\beta}^\circ=\ker(\varphi_{x_\beta})$.
\end{proof}

We now prove
Theorems \ref{thm:E7_ROD}, \ref{thm:E7_orbit_parameterization}
using Propositions
\ref{prop:E7_neutral_element},
\ref{prop:E7_quadratic_element},
\ref{prop:E7_cubic_element}
and Corollary \ref{cor:fiber_structure}.
\begin{proof}
We already know the map $\Psi$ is surjective
and we consider the fiber structure.
Let $x\in V^\sst(k)$ be either one of
$w$, $x_\alpha$ or $x_\beta$ in
\eqref{eq:E7_w}, \eqref{eq:E7_xalpha} or \eqref{eq:E7_xbeta}.
For such $x$, we proved in Propositions
\ref{prop:E7_neutral_element},
\ref{prop:E7_quadratic_element} and
\ref{prop:E7_cubic_element} that
$G_x^\circ$ satisfies the exact sequence
\[
1
\longrightarrow		G_x^\circ
\longrightarrow		K^x
\overset{\varphi_x}{\longrightarrow}	\esR_{L_x/k}({\mathbb G}_m)
\longrightarrow		1
\]
where $K^x\cong {\mathbb G}_m\times \esR_{L_x/k}(\esB_{L_x}^\times)$.
By Proposition \ref{prop:cohomology_compute},
$\rh^1(k,K^x)=\{1\}$.
Hence by Theorem \ref{thm:rationalorbits_exact} we have
\begin{equation}\label{eq:parameter_set}
\rh^1(k,G_x^\circ)\cong \varphi_x(K^x(k))\backslash L_x^\times.
\end{equation}
Since $s$ is a contravariant functor,
\begin{equation}\label{eq:auto_group}
\aut(Z_x)(k)=\aut_k(Z_x)\longrightarrow
\aut_k(L_x),
\qquad
 \psi\longmapsto s^\ast(\psi^{-1})
\end{equation}
gives the isomorphism of the groups.
We later see their actions on respective sets are compatible.
With the compatibility,
by Corollary \ref{cor:fiber_structure} we have the parameterization.

We give the explicit correspondence of this parameterization
and determine the structure of stabilizers of the elements.
First consider the case $x=w$, $L=L_w=k\times k\times k$.
Take $p=(p_1,p_2,p_3)\in L^\times$.
Let $g=\sfd(b_1,b_2,b_3,1)\in K^w(\ks)$ be in $\varphi_w^{-1}(p)$,
i.e., $\cN(b_i)=p_i$ for $i=1,2,3$.
Then the class of $(g^{-1}g^\sigma)_\sigma$ corresponds to the class of $p$.
Hence by Corollary \ref{cor:fiber_structure}, the class of $p\in L^\times$
corresponds to $gw$, which equals to the element of
Theorem \ref{thm:E7_orbit_parameterization} (1).
Also by the form of $gw$ elements of $G_w^\circ$ stabilize $gw$.
Comparing the dimensions, we have $G_{gw}^\circ=G_w^\circ$.

Next consider the case $x=x_\beta$, $L=L_{x}$ is
a cubic Galois extension of $k$.
Let $\theta\in \gal(L_x/k)$ be the element satisfying
$\beta_2=\beta_1^\theta$, $\beta_3=\beta_2^{\theta}$.
Let $\delta\in L^\times$.
We use the notation $\delta_i$ and $\Delta_i$
as in Theorem \ref{thm:E7_orbit_parameterization} (3).
Take $b\in (\esB_L\otimes\ks)^\times$ such that $\cN(b)=\delta$
and put $g=g_\beta\sfd(b,b^\theta,b^{\theta^2},1)g_\beta^{-1}$.
Then $\varphi_x(g)=\delta$ and so
the class of $(g^{-1}g^\sigma)_\sigma$ corresponds to
the class of $\delta$.
Hence by Corollary \ref{cor:fiber_structure}, 
the class of $\delta\in L^\times$
corresponds to $gx$. Since $x=g_\beta w$,
\[
gx=g_\beta\sfd(b,b^\theta,b^{\theta^2},1)w
=g_\beta
\left(\thth 0000{\delta^\theta}000{-\delta^{\theta^2}},
\thth \delta000{-\delta^\theta}0000
\right)
\]
and thus equals to the element of
Theorem \ref{thm:E7_orbit_parameterization} (3).
Since elements of $G_{x}^\circ$ stabilize $gx$,
we have $G_{gx}^\circ=G_x^\circ$.
The remaining cases are proved similarly and we omit the details here.

We finally give the compatibility of actions of groups
\eqref{eq:auto_group} on sets \eqref{eq:parameter_set}.
Since the proofs are similar we consider the case $L=L_{x}$ is
a cubic Galois extension of $k$.
We use the notation above.
We put
\[
\tilde\theta=
\left(\thth 001100010,\twtw 01{-1}{-1}\right)\in G(k).
\qquad
\]
Then $g_\beta^\theta=g_\beta\tilde\theta$
and $\tilde\theta\in G_w$.
Hence if we let
$h=g_\beta\tilde\theta g_\beta^{-1}$,
then we have $h\in G_x(k)$.
Recall that we define $c_x:G_x\rightarrow \aut(Z_x)$
in \eqref{eq:c_w}.
By definition,
$Z_x=\{(\beta_1:-1),(\beta_2:-1),(\beta_3:-1)\}\subset \mathbb P^1$.
We see by computation that $c_x(h):Z_x\rightarrow Z_x$ is given by
\[
(\beta_1:-1)\mapsto(\beta_3:-1),
\quad
(\beta_2:-1)\mapsto(\beta_1:-1),
\quad
(\beta_3:-1)\mapsto(\beta_2:-1).
\]
Hence under the isomorphism in \eqref{eq:auto_group},
$\theta\in\aut(L_x)$ corresponds to $c_x(h)\in\aut(Z_x)(k)$.
By the definition of the action,
$(g^{-1}g^\sigma)_\sigma^{c_x(h)}
=(h^{-1}g^{-1}g^\sigma h^\sigma)_\sigma
=((h^{-1}gh)^{-1} (h^{-1}gh)^\sigma)_\sigma$.
Since
\[
h^{-1}gh
=	g_\beta
	\tilde\theta^{-1}\sfd(b,b^\theta,b^{\theta^2},1)\tilde\theta
	g_\beta^{-1}
=	g_\beta
	\sfd(b^\theta,b^{\theta^2},b,1)
	g_\beta^{-1},
\]
the image of $\delta^\theta$ under the bijection \eqref{eq:parameter_set}
is the cohomology class of the cocycle $(g^{-1}g^\sigma)_\sigma^{c_x(h)}$.
This shows the compatibility of the actions.
\end{proof}

\subsection{Decomposition of the zeta function}
We conclude this section with a brief discussion
of the global zeta functions of the prehomogeneous vector space $(G,V)$.
Let $k$ be a number field and $\A$ its adele ring.
Let $(G,V)$ be either the $D_5$ or $E_7$ type.
We put $\cA_2'=\{L\in\cA_2\mid \cA_2\neq k\times k\}$
and $\cA_3'=\cA_3$.
Let $\esL=\{x\in V^\sst(k)\mid k(x)\in\cA_n'\}$.
The global zeta function for $\Phi=\prod_v \Phi_v$
in the Schwartz-Bruhat space $\sS(V(\A))$ is defined by
\[
Z(\Phi,s)=\int_{G(\A)/G(k)}|\chi(g)|_\A^{2s}\sum_{x\in\esL}\Phi(gx)dg.
\]
It is proved by H. Saito \cite{hsaitoc}
that this integral converges absolutely for
sufficiently large $\re(s)$.
For $L\in \cA_n$, let $V_{k,L}=\{x\in V^\sst(k)\mid k(x)=L\}$.
Then since $\esL=\coprod_{L\in\cA_n'}V_{k,L}$ we have
\begin{align*}
Z(\Phi,s)
&	=\sum_{L\in\cA_n'}Z(\Phi,s,L),\\
Z(\Phi,s,L)
&	=\int_{G(\A)/G(k)}|\chi(g)|_\A^{2s}\sum_{x\in V_{k,L}}\Phi(gx)dg.
\end{align*}
If $(G,V)$ is of the $D_5$ type, $V_{k,L}$ is a single $G(k)$-orbit.
Also since the first Galois cohomology of $G_x^\circ$ vanishes
for any local field of $k$,
we have $G(\A)/G_x^\circ(\A)=(G/G_x^\circ)(\A)$.
These combined with the usual unfolding method show that
$Z(\Phi,s,L)$ has an Euler product and each Euler factor
is what is called the local zeta function.
H. Saito \cite{hsaitob} proved that
even if $G(\A)/G_x^\circ(\A)\neq(G/G_x^\circ)(\A)$ for some $x\in V^\sst(k)$,
the global zeta function $Z(\Phi,s)$
has a formula expressed via local zeta functions.
Let $(G,V)$ be of the $E_7$ type.
If $k$ has a real place, then $G(\A)/G_x^\circ(\A)\neq(G/G_x^\circ)(\A)$
for certain $x\in V^\sst(k)$.
But the decomposition $\esL=\coprod_{L\in\cA_n'}V_{k,L}$
is according to $\h^1(k,G_x/G_x^\circ)$ for an element $x\in V^\sst(k)$.
Hence his result shows that
$Z(\Phi,s,L)$ is a finite number of functions
where each of them has an Euler product,
and each Euler factor is a finite sum of local zeta functions.

\section{Parameterizations of reducible algebras for cubic cases}
\label{sec:reducible}
In \cite{wryu}, Wright-Yukie showed that the
semi-stable rational orbits $G(k)\backslash V^\sst(k)$
of the following three prehomogeneous vector spaces
\begin{center}
\begin{tabular}{llll}
(1) & $G=\gl_1\times\gl_2$, & $V=\sym^3k^2$, & ($G_2$ type)\\
(2) & $G=\gl_3\times\gl_3\times\gl_2$, & $V=k^3\otimes k^3\otimes k^2$,
& ($E_6$ type)\\
(3) & $G=\gl_6\times\gl_2$, & $V=\wedge^2k^6\otimes k^2$,
& ($E_7$ type)
\end{tabular}
\end{center}
parameterize separable cubic algebras $\cA_3$ of $k$.
Let $\cA_3=\cE_3\amalg \cE_{1,2}\amalg \cE_{1,1,1}$, where
\begin{align*}
\cE_3&=\{L\mid \text{$L/k$ is a cubic field extension}\},\\
\cE_{1,2}&=\{L=k\times F \mid \text{$F/k$ is a quadratic field extension}\},\\
\cE_{1,1,1}&=\{k\times k\times k\}.
\end{align*}
From number theoretical interests,
we often would like to know about the subset $\cE_3$ rather than
the whole set $\cA_3$ itself.
Let $\br_2\subset\gl_2$ denote the set of lower triangular matrices.
For the classical case (1), Shintani \cite{shintanib}
introduced the representation $(\br_2,\sym^2k^2)$
to parameterize $\cE_{1,2}\amalg \cE_{1,1,1}$,
and found a good information on $\cE_3$.
Shintani considered $\Z$-orbits and his method
was recently extended to over a Dedekind domain by the author \cite{ddca}.
In this section, we find analogous representation
parameterizing $\cA_3\setminus\cE_3$ for (2), (3) and
their $k$-forms over fields.
The study of integral orbits will be important.
Also similar problems should arise for the quartic case
$(\gl_3\times\gl_2,\sym^2 k^3\otimes k^2)$
and quintic case $(\gl_5\times\gl_4, \wedge^2k^5\otimes k^4)$.
These would be answered in the future.

For a quaternion algebra $\esB$ of $k$ we keep the notation
in Section \ref{sec:rationalorbitsD5E7}.
Let $F$ be a quadratic field extension of $k$.
We denote by $\sigma$ the generator of $\gal(F/k)$.
Let ${\rm N}_{F/k}$ be the norm map.
Let $\rh_3(F)=\{x\in {\rm M}_3(F)\mid {}^tx=x^\sigma\}$.
The group $\gl_3(F)$ acts on this space
via $(g,x)\mapsto gx{g^\sigma}$, for $g\in\gl_3(F)$, $x\in \rh_3(F)$.

We define and discuss the three spaces in parallel;
they will be distinguished as Cases (a), (b) and (c).
Let
\begin{center}
\begin{tabular}{llll}
$G_1=\gl_3(k)\times\gl_3(k)$, &&$Y={\rm M}_3(k)$	& in Case (a),\\
$G_1=\gl_3(F)$,		&&	$Y=\rh_3(F)$		& in Case (b),\\
$G_1=\gl_3(\esB)$,	&&	$Y=\rh_3(\esB)$		& in Case (c),\\
\end{tabular}
\end{center}
and put
\[
G=G_1\times\gl_2(k),\qquad V=Y\oplus Y.
\]
(We already defined and considered Case (c)
 in Section \ref{sec:rationalorbitsD5E7}.
For later conveniences we do not admit $F$ to be $k\times k$
and treat Cases (a) and (b) separately.)
We regard $G$ as an algebraic group over $k$
and $V$ as a variety over $k$.
We write elements of $G$ as $g=(g_{11},g_{12},g_2)\in G$ in Case (a)
and $g=(g_1,g_2)$ in Cases (b) and (c).
We identify elements $x=(x_1,x_2)\in V$
with the $3$ by $3$ matrices $x(v)=x_1v_1+x_2v_2$ 
of linear forms in the variables $v_1$ and $v_2$,
which we collect into the row vector $v=(v_1,v_2)$.
With this identification, we define a rational action of $G$ on $V$ via
\[
(g\cdot x)(v)
=\begin{cases}
g_{11}x(vg_2)^{t}\!g_{12}	&	\text{in Case (a)},\\
g_1x(vg_2)^{t}\!g_1^\sigma	&	\text{in Case (b)},\\
g_1x(vg_2)g_1^\iota		&	\text{in Case (c)},\\
\end{cases}
\]
Let $F_x(v)=\det(x(v))$ in Cases (a), (b)
and $F_x(v)=\pfaff_3(x(v))$ in Case (c).
Let $P(x)$ be the discriminant of the binary cubic form $F_x(v)$.
Then $P\in k[V]$ is a homogeneous polynomial of degree 12
and $P(g\cdot x)=\chi(g)^2P(x)$ for $g\in G$, $x\in V$ where
\[
\chi(g)=
\begin{cases}
\det(g_1)^2\det(g_2)^2\det(g_3)^3	&	\text{in Case (a)},\\
{\mathrm N}_{F/k}(\det(g_1))^2\det(g_2)^3	&	\text{in Case (b)},\\
\cN(g_1)^2\det(g_2)^3			&	\text{in Case (c)}.\\
\end{cases}
\]
Let $V^\sst=\{x\in V\mid P(x)\neq0\}$.
\begin{defn}
For $x\in V^\sst(k)$, we define
\[
Z_x=\mathrm{Proj}\, k[v_1,v_2]/(F_x(v)),\quad
L_x=\Gamma(Z_x,\co_{Z_x}).
\]
We define
\begin{align*}
V^{(3)}&=\{x\in V^\sst(k)\mid L_x\in\cE_3\},\\
V^{(2)}&=\{x\in V^\sst(k)\mid L_x\in\cE_{1,2}\},\\
V^{(1)}&=\{x\in V^\sst(k)\mid L_x\in\cE_{1,1,1}\}
\end{align*}
\end{defn}
The algebra $L_x\in\cA_3$ depends only on $G(k)$-orbit of $x$.
Hence each $V^{(i)}$ is a $G(k)$-invariant subset of $V^\sst(k)$.

For $A=k,F,\esB$, we put
\[
{\mathrm P}_{1,2}(A)=\left\{\thth *00******\in\gl_3(A)\right\}.
\]
Let $\br_2\subset \gl_2$ be the set of lower triangular matrices.
We define
\[
P=
\begin{cases}
{\mathrm P}_{1,2}(k)\times {\mathrm P}_{1,2}(k)\times \br_2(k)
&	\text{in Case (a)},\\
{\mathrm P}_{1,2}(F)\times \br_2(k)
&	\text{in Case (b)},\\
{\mathrm P}_{1,2}(\esB)\times \br_2(k)
&	\text{in Case (c)}.\\
\end{cases}
\]
We regard $P$ as an algebraic group over $k$,
then is a parabolic subgroup of $G$.
For $Y={\mathrm M}_3(k)$, $\rh_3(F)$ or $\rh_3(\esB)$, we put
\[
Z=\left\{\thth 0000**0**\in Y\right\}
\]
and define $W=Z\oplus Y\subset V$.
Then $W$ is a $P$-invariant subspace
and $(P,W)$ is a prehomogeneous vector space
with $W^\sst=W\cap V^\sst$ a single orbit over $\overline k$.
We put $W^{(i)}=W^\sst(k)\cap V^{(i)}$ for $i=1,2$.
We note that for $w\in W^\sst$,
$F_w(v)$ has a linear factor $v_2$ and hence 
$W^\sst(k)\cap V^{(3)}=\emptyset$.

We first consider the relation between $V^{(2)}$ and $W^{(2)}$.
\begin{prop}
We have
\[
V^{(2)}=G(k)\times_{P(k)}W^{(2)}.
\]
\end{prop}
\begin{proof}
In Case (c), we show in (2) of Theorem \ref{thm:E7_orbit_parameterization}
that any orbit in $V^{(2)}$ has an element in $W^\sst(k)$
and hence $V^{(2)}=G(k)\cdot W^{(2)}$.
We did not consider Cases (a) and (b) in Section
\ref{sec:rationalorbitsD5E7} but exactly the same arguments works and
we can see any orbit in $V^{(2)}$ has an element in $W^\sst(k)$.
Hence again we have $V^{(2)}=G(k)\cdot W^{(2)}$.
The rest of arguments are also parallel for Cases (a), (b) and (c),
and we write proofs for Case (c).

It is obvious that $P(k)$ fixes the set $W^{(2)}$.
Let $w\in W^{(2)}$, $g=(g_1,g_2)\in G(k)$ and $x=g\cdot w\in W^{(2)}$.
We will show $g\in P(k)$.
We first show $g_2\in\br_2(k)$.
Since $x=g\cdot w$, $P_x(v)$ and $P_w(vg_2)$ coincide
up to a constant multiple.
On the other side since $w\in W^{(2)}$,
$P_w(v)$ is a product of $v_2$ and
an irreducible quadratic polynomial in $v$, and $P_x(v)$ is also.
Hence $g_2\in\gl_2(k)$ must fix $(1:0)\in\mathbb P^1(k)$ and
consequently $g_2\in\br_2(k)$.

Let
\[
w=
\left(
	\left(
	\begin{array}{c|c}0&0\\\hline0&w_{122}\\\end{array}
	\right),
	\left(
	\begin{array}{c|c}w_{211}&w_{212}\\\hline w_{221}&w_{222}\\\end{array}
	\right)
\right),
\quad
w_{122},w_{222}\in\rh_2(\esB), \ w_{211}\in k,
\]
be the block representation of $w=(w_1,w_2)$.
Then the coefficients of $v_1^2v_2$ of $F_w(v)$
is $w_{211}\pfaff_2(w_{122})$ and this does not vanish since $w\in W^\sst(k)$.
Hence $w_{122}$ is invertible i.e., $w_{122}\in\gl_2(\esB)$.
We consider the similar block representation for $x$, and we have
$x_{122}\in\gl_2(\esB)$.
Let $g_2=\stwtw a0bc$, $p=ag_1$, $q=(g_1^\iota)^{-1}$.
Then since $x_1=ag_1w_1g_1^\iota$ we have $x_1q=pw_1$.
Write $p,q\in\gl_3(\esB)$ in the block form
\[
p=	\left(
	\begin{array}{c|c}p_1&p_2\\\hline p_3&p_4\\\end{array}
	\right),
\quad
q=	\left(
	\begin{array}{c|c}q_1&q_2\\\hline q_3&q_4\\\end{array}
	\right),
\qquad
(p_1,q_1\in\esB,\ p_4,q_4\in{\mathrm M}_2(\esB).)
\]
Then $p_2x_{122}=(0\ 0)$ and $w_{122}q_3 ={}^t(0\ 0)$.
Since $x_{122}, w_{122}\in\gl_2(\esB)$, we have
$p_2=(0\ 0)$, $q_3={}^t(0\ 0)$.
Thus we have $g=(g_1,g_2)\in P(k)$ and this finishes the proof.
\end{proof}

We next consider the relation between $V^{(1)}$ and $W^{(1)}$.
For $A=k,F,\esB$, we put
\[
{\mathrm T}_3(A)=\left\{\thth *000*000*\in\gl_3(A)\right\}.
\]
Let $\mathrm Z_2\subset \gl_2$ be the center of $\gl_2$.
We define
\[
H^\circ=
\begin{cases}
{\mathrm T}_{3}(k)\times {\mathrm T}_{3}(k)\times {\mathrm Z}_2(k)
&	\text{in Case (a)},\\
{\mathrm T}_{3}(F)\times {\mathrm Z}_2(k)
&	\text{in Case (b)},\\
{\mathrm T}_{3}(\esB)\times {\mathrm Z}_2(k)
&	\text{in Case (c)}.\\
\end{cases}
\]
We regard $H^\circ$ as an algebraic group over $k$.
We put
\[
\theta_1=\thth 001100010,
\quad
\theta_2=\twtw 01{-1}{-1},
\quad
\eta_1=\thth 100001010,
\quad
\eta_2=\twtw {-1}011
\]
and define
\[
\theta=
\begin{cases}
(\theta_1,\theta_1,\theta_2)
&	\text{in Case (a)},\\
(\theta_1,\theta_2)
&	\text{in Case (b)},\\
(\theta_1,\theta_2)
&	\text{in Case (c)},\\
\end{cases}
\qquad
\eta=
\begin{cases}
(\eta_1,\eta_1,\eta_2)
&	\text{in Case (a)},\\
(\eta_1,\eta_2)
&	\text{in Case (b)},\\
(\eta_1,\eta_2)
&	\text{in Case (c)}.\\
\end{cases}
\]
Let $\gA_3$, $C_2$ be respectively the group generated by
$\theta$ and $\eta$, and $\gS_3$ the group generated by $\gA_3$ and $C_2$.
Then $\gS_3$ is isomorphic to the symmetric group of degree 3,
$\gA_3$ to the alternating subgroup,
and $\gS_3=\gA_3\rtimes C_2$.
It is easy to see that $\gS_3$ normalizes $H^\circ$.
We put $H=H^\circ\rtimes\gS_3$. Then $H^\circ$
is the identity component of $H$.
We put
\[
U=\left\{\left(\thth 0000b000{-c}, \thth a000{-b}0000\right)\in V\ 
\vrule\ a,b,c\in k\right\}.
\]
Then $U$ is a $H$-invariant subspace. We put $U^\sst=U\cap V^\sst$.
For $u\in U$, $F_u$ coincides with $v_1v_2(v_1-v_2)$ up to constant multiple.
Hence $U^\sst(k)\subset W^{(1)}$.
\begin{prop}
We have
\begin{align*}
V^{(1)}&=G(k)\times_{H(k)}U^{\sst}(k),\\
W^{(1)}&=P(k)\times_{H(k)\cap P(k)}U^{\sst}(k).
\end{align*}
\end{prop}
\begin{proof}
Since the proofs are similar, we here write a proof for Case (c).
Let us consider the first formula. 
We show in (1) of Theorem \ref{thm:E7_orbit_parameterization}
that any orbit in $V^{(1)}$ has an element in $U^\sst(k)$
and hence $V^{(1)}=G(k)\cdot U^\sst(k)$.
Let $x\in U^\sst(k)$, $g=(g_1,g_2)\in G(k)$ and assume $y=gx\in U^\sst(k)$.
Since the set of roots of $F_x(v)$ and $F_y(v)$ in $\mathbb P^1$ are both
$\{(0:1), (1:0), (1:1)\}$, $g_2\in\gl_2(k)$ fixes this set.
Since the $\gl_2$-part of $\gS_3$ coincides with
the permutation group of this set, there exists
$\sigma=(\sigma_1,\sigma_2)\in\gS_3$
such that $\sigma_2g_2\in\gl_2(k)$ is a scalar matrix.
Hence $\sigma gx\in U^\sst(k)$ and the $\gl_2$-component of $\sigma g$
is in $\mathrm Z_2(k)$.
Now it is easy to see that $G_1$-component of $\sigma g$
is in $T_3(\esB)$. Hence $\sigma g\in H^\circ(k)$ and we have $g\in H(k)$.
This proves the first formula.

For the second formula, it is enough to show $W^{(1)}=P(k)\cdot U^\sst(k)$.
Let $w\in W^{(1)}$ and
\[
w=
\left(
	\left(
	\begin{array}{c|c}0&0\\\hline0&w_{122}\\\end{array}
	\right),
	\left(
	\begin{array}{c|c}w_{211}&w_{212}\\\hline w_{221}&w_{222}\\\end{array}
	\right)
\right),
\quad
w_{122},w_{222}\in\rh_2(\esB), \ w_{211}\in k,
\]
be the block representation of $w=(w_1,w_2)$.
We can prove $w\in P(k)\cdot U^\sst(k)$ by the following steps.
Let $P_{1,2}(\esB)=L_{1,2}(\esB)\ltimes U_{1,2}(\esB)$
be the Levi decomposition.
\begin{enumerate}
\item Since $w_{211}\in \mk$, by multiplying an element
of $U_{1,2}(\esB)$ if necessary, we may assume $w_{212}=w_{221}=0$.
\item Because of the result of $D_5$-type, by multiplying an element
of $L_{1,2}(\esB)$ if necessary, we may 
further assume $w_{122}$ and $w_{222}$
are both diagonal.
\item Moreover, by multiplying an element of $B_2(k)$,
we may assume the set of roots of $F_w(v)$ is $\{(0:1), (1:0), (1:1)\}$.
\end{enumerate}
Under these conditions,
either $w$ or $\eta w$
must be an element of $U^\sst(k)$.
This finishes the proof.
\end{proof}

\end{document}